\documentclass[twoside]{article}
\usepackage{color,times,amsmath,amsfonts,latexsym,epsfig,graphicx,epsf,colordvi}
\usepackage{url,hyperref}
\usepackage[english]{babel}

\usepackage{amssymb}
\usepackage{comment}
\usepackage{url}
\usepackage{fancyhdr}
\usepackage{mathtools}

\def\tcb{\textcolor[rgb]{0,0,1}}
\def\tcr{\textcolor[rgb]{1,0,0}}

\newcommand{\CC}{\mathbb {C}}
\newcommand{\RR}{\mathbb {R}}

\newcommand{\la}{\lambda }

\newcommand{\et}{\eta }
\newcommand{\bp}{\begin{pmat}}
\newcommand{\ep}{\end{pmat}}

\makeatletter
\def\adots{\mathinner{\mkern1mu\raise\p@\vbox{\kern7\p@\hbox{.}}\mkern2mu\raise4\p@\hbox{.}\mkern2mu\raise7\p@\hbox{.}\mkern1mu}}
\makeatother

\title{ \vspace*{-8mm}{ \bf Time-Varying Matrix Eigenanalyses via Zhang Neural Networks and look-Ahead Finite Difference Equations \vspace*{-4mm} }}

\author{Frank Uhlig\\
Department of Mathematics and Statistics,\\
Auburn University, AL 36849-5310, USA \ ({\tt uhligfd@auburn.edu})\\[2mm]
Yunong Zhang\\
School of Information Science and Technology,\\ 
Sun Yat-sen University (SYSU),\\ 
Guangzhou 510006, China \ ({\tt zhynong@mail.sysu.edu.cn})\\[2mm]}

\begin{document}
\date{  }
\maketitle

\vspace*{-2mm}

\pagestyle{fancy}
\fancyhead{}
\fancyhf{} % sets both header and footer to nothing
\renewcommand{\headrulewidth}{0pt}

\fancyhead[CE]{Frank Uhlig and Yunong Zhang}
\fancyhead[CO]{Time-varying Matrix Eigenvalues via ZNN and Finite Difference Equations} 
\fancyhead[RO]{\thepage}
\fancyhead[LE]{\thepage}
\thispagestyle{empty}

\vspace*{-6mm}
{\normalsize
\noindent 
{\bf Abstract : }
\noindent
This paper adapts look-ahead and backward finite difference formulas to compute future eigenvectors and eigenvalues of piecewise smooth time-varying symmetric matrix flows $A(t)$. It is based on the Zhang Neural Network (ZNN) model for time-varying  problems and uses the associated error function $E(t) = A(t)V(t) - V(t) D(t)$ or $e_i(t) = A(t)v_i(t) -\la_i(t)v_i(t)$ with the Zhang design stipulation that $\dot E(t) = - \eta E(t)$ or $\dot e_i(t) = - \eta e_i(t)$ with $\eta > 0$ so that $E(t)$ and $e(t)$ decrease exponentially over time. This leads to a discrete-time differential equation of the form $P(t_k) \dot z(t_k) = q(t_k)$ for the eigendata vector $z(t_k)$ of $A(t_k)$. Convergent look-ahead finite difference formulas of varying error orders then allow us to express $z(t_{k+1})$ in terms of earlier $A$ and $z$ data. Numerical tests, comparisons and open questions complete the paper.\\[2mm]
{\bf Subject Classifications :} \ 65H17, 65L12, 65F15,  65Q10, 92B20 \\[2mm]
{\bf Key Words :}  matrix eigenvalues, matrix computations, symmetric matrix flows, Zhang neural network, time-varying equations, time-varying computations, 1-step ahead discretization formula,  error function,  numerical analysis, numerical linear algebra}

\section{\Large Introduction}

Many numerical methods in most every branch of applied mathematics  employ matrices, vectors and spatial or time-dependent models to solve specific problems. Static-time and time-varying problems sometimes behave very differently as far as their numerics are concerned. Therefore time-invariant and time-varying problems may require different approaches when they deal with different  challenges and their own inherent computational limitations.

Time-varying numerical matrix matrix and ZNN methods are a relatively new subject and differ greatly from our classical numerical canon and textbooks. The very words \emph{neural networks} have multiple uses in Numerics and simply  refer to the natural propagation of impulses that are passed along nervous systems. ZNN  methods differ substantially from every other so called neural network, see e.g. \cite{WCW16}, and from any other standard ODE based method, see e.g. \cite{LM}, or from decomposition methods, see e.g. \cite{DE99}.\\[1mm] 
ZNN methods have  two main and differing ingredients: the first is the {\em error equation} and the {\em ZNN stipulation} that the error equation is to decay exponentially fast to zero over time. Once this new,  quite unusual and nowhere else used error function ODE has been discretized in ZNN, the second  novel part of ZNN is its subsequent reliance on {\em 1-step ahead convergent finite difference equations}, rather than on standard ODE IVP solvers. Such look-ahead difference schemes  have never occurred in any of  our finite difference uses or the literature. The first ones were  in fact constructed by hand on scratch paper within the last couple of years, see \cite{LMUZb2018} and \cite{QZY2019}.  These first tries achieved rather low truncation error orders. A formal process to obtain higher error order 1-step ahead convergent difference schemes appears in \cite{U2019}.\\[1mm]
There are several hundred papers with time-varying Zhang type methods in  engineering journals, but relatively few, not even a handful, in  the Numerical Linear Algebra literature.\\[1mm]
 For over a dozen years now, a special class of dynamic methods has been built on the idea of Yunong  Zhang and Jun Wang \cite{ZW01} from 2001. These are so called \emph{Zhang dynamics} (ZD), or \emph{zeroing dynamics}; see \cite{ZKXY2010,ZLYL2012,GNY2016,XL2016} and \emph{Zhang neural networks} (ZNN), or \emph{zeroing neural networks}; see \cite{ZY2011,JLH2018,XLLZD2018}. ZD methods are specifically designed for time-varying problems and have proven most efficient there. They  use   1st-order time derivatives and have been applied successfully to  -- for example -- solve  time-varying Sylvester equations  \cite{ZJW2002,LL2014,LCL2013}, to  find time-varying matrix inverses \cite{ZG2005,JLH2018,C2013,LMUZb2018} (see also \cite{CY2016}), and to optimize  time-varying matrix minimization problems \cite{ZYLHW2017,LMUZa2018}, all in real-time. These algorithms generally use  both, 1-step ahead and backward differentiation formulas and  run in discrete time with high  accuracy. From  given time-varying matrix or vector valued problem, all ZNN methods form a problem specific  error function $e(t)$  that  is  stipulated  to decrease exponentially to zero, both globally and asymptotically \cite{ZLYL2012,ZY2011,ZYGZ2011} by the ODE $ \dot e(t) = -\et e(t)$ for a decay constant $\et > 0$. This  discretized time-varying error matrix ODE problem is then solved successfully through finite difference formulas and simple recursive vector additions without the use of any  standard numerical matrix factorizations, local iterative methods, or ODE solvers, except  possibly for  computing the necessary starting values.\\[1mm]
  The ZNN method is a totally new branch of Numerical Analysis. It is essential for robot control, self-driving vehicles, and autonomous airplanes et cetera where its ability to predict future systems states accurately is an essential ingredient for engineering success. ZNN's numerical findings and open challenges are many and most of these have not even been recognized or even named in their own right. \\[1mm]
In this paper we extend  the ZD method for the recently studied time-varying matrix eigenvalue  problem \cite{ZYLUH2018} that dealt with convergence and robustness issues of the DE based ZNN model solution. Here we propose a new symmetric discrete-time ZD model that computes all  eigenvalues and  eigenvectors of  time-varying real symmetric matrix flows $A(t)$, be they repeating or not. Our model  uses   convergent look-ahead and standard backward finite difference formulas instead of an ODE solver as was done  in \cite{ZYLUH2018} for example.\\ 
Previous efforts to find the eigenstructure of time-varying matrix flows, mostly  via DE solvers and path following continuations go back at least  to 1991 \cite{BBMN91}  where the smoothness of the eigendata and SVD data of a time-dependent matrix flow was explored in detail. A few years later, Mailybeav \cite{M06} studied  problems with near identical eigenvalues for time-varying general matrix computations. Matrix flows with coalescing eigenvalues pose {\emph{crucial problems}} for these methods. These problems were observed and further studied by Dieci and Eirola in 1999 \cite{DE99} and more recently by  Sirkovi\'c and Kressner   \cite{SK16} in 2016 and they likewise occur for a related parameter-varying matrix eigen problem of Loisel and Maxwell in 2018 \cite{LM}. Our ZNN eigenvalue algorithm is impervious to these restrictions as it handles repeated eigenvalues in symmetric matrix flows $A(t)$ without any problems, see Figures 1, 5, and 8.

The new ZD model of this paper is a discrete dynamical system with the potential for practical real-time and on-chip implementation and computer simulations. We include results of computer simulations and numerical experiments that illustrate the usefulness and efficiency of our real-time discrete ZD matrix eigenvalue algorithm, both for  smooth data inputs and also for piecewise smooth time-varying matrix flows $A(t)$ that might occur naturally when sensors fail or data lines get disrupted in the field.

\section{\Large A Zhang Neural Network Formulation for the Time-Varying Matrix Eigenvalue Problem}

Computing the eigenvalues and eigenvectors of matrices with known fixed entries , i.e., of static matrices, has presented mathematicians with difficult problem from the early 1800s on. In the early 20th century, engineers such as Krylov thought of vector iteration and built iterative methods for finding matrix eigenvalues. Iterative matrix methods increased in sophistication  with the advent of computers in the 1950s and are now very successful to solve static matrix problems for a multitude of  dense or large structured and sparse matrices. For dense matrices the eigenvalue problem was  essentially left unsolved for 150 years until John Francis and Vera Kublanovskaja independently created the QR algorithm around 1960.   \\[1mm]
While it is nowadays quite easy to solve the eigenvalue problem $Ax\ = \la x$ for  fixed entry  matrices $A$ and find their eigenvectors $x_i$ and eigenvalues $\la_i$, what happens if the entries of $A$ vary over time? Naively, one could solve the associated time-varying  eigenvalue equation $A(t)x(t) = \la(t) x(t)$ for multiple times $t_k$ and take the computed eigendata output for $A(t_k)$ as an approximate  solution of the time-varying eigen problem at time $t_{k+1}$. This might work if $A(t)$ is explicitly known for all $t_0 \leq t\leq t_{end}$ and if we can wait for the end of the computations for $A(t_k)$ to see the results for the past time $t_k$. 
 But if we need to know the eigenstructure of $A(t_{k+1})$ from the previous behavior of $A(t)$ before time $t_{k+1}$ has arrived, such as for robot optimization processes etc,  the naive  method would compute  and deliver the 'exact' solution for an instance in the past, one  that may have little or no relation to the current time situation. And besides, how would the delayed solution maneuver discontinuities in $A(t)$ in real-time or when large entry value swings occur in $A(t)$? Therefore to solve time-varying matrix eigenvalue problems reliably in real-time we need  to  learn how to \emph{predict the future value} of the eigenvalues and eigenvectors of $A(t)$ for the time instance $t_{k+1}$ solely from earlier and past  $A(t_{j})$ data  instances $t_j$ with $t_j < t_{k+1}$. This  challenge is rather new and has not been much explored in numerical analysis.\\[1mm]
 Our best static matrix eigenvalue methods use orthogonal matrix factorizations and these are backward stable. Backward stability gives us the exact solution to an adjacent problem whose distance from the given problem is bounded by the eigenvalue condition numbers of the  given matrix.
For time-varying matrix problems, we apparently need to relax  on backward stability and instead work on accurate forward predictions. Here is how this can be done by using Zhang Neural Network ideas for  time-varying matrix eigenvalue problems.\\[0mm]

{\bf \large Statement of the Problem and the Zhang Neural Network Error Equation}\\[2mm]
For a real symmetric flow of matrices $A(t)\in \RR^{n,n}$  and $0 \leq t \leq  t_f$, we consider the problem of finding  non-singular real matrices $V(t)\in \RR^{n,n}$ and  real diagonal matrices $D(t)\in \RR^{n,n}$ so that 
\begin{equation}\label{EVequat}
A(t) V(t) = V(t)D(t) \ \text{ for all } t \ .
\end{equation}
Since we assume $A(t)$ to be  real symmetric at all times $t$, such matrices $V(t)$ and $D(t)$ will  exist for all $t$. Our aim is to  find them accurately and predictively in real-time. To solve (\ref{EVequat}), the Zhang Neural Network approach  looks at the time-varying homogeneous error equation which  has the form
\begin{equation}\label{errequat}
e(t) = A(t)V(t) - V(t) D(t)  \stackrel{!}{=} O_n \in \RR^{n,n} \ .
\end{equation}
Next, the ZNN approach  stipulates exemplary behavior for $e(t)$ by asking for exponential decay  of $e(t)$ as a function of time, or
\begin{equation}\label{errde}
\dot e(t) = - \eta e(t)
\end{equation}
with a decay constant $\eta > 0$. Equation (\ref{errde}) can be written out explicitly as
\begin{equation}\label{errdeexpl}
\dot e(t) = \dot A(t) V(t) + A(t)\dot V(t) - \dot V(t) D(t) - V(t) \dot D(t) =  - \eta A(t)V(t) + \eta V(t)D(t) = - \eta e(t) , \text{ or} 
\end{equation}
\begin{equation}\label{errdeexpl2}
  A(t) \dot V(t) - \dot V(t) D(t) - V(t)  \dot D(t) =  - \eta A(t)V(t) + \eta V(t)D(t)  - \dot A(t) V(t)
 \end{equation}
where in (\ref{errdeexpl2}) we have gathered all terms with the derivatives $\dot V(t)$ of the unknown eigenvector matrix $V(t)$ and $\dot D(t)$ of the eigenvalue matrix $D(t)$ on the left hand side.\\
When we specify equation (\ref{errdeexpl2})  for one eigenvalue/eigenvector pair $x_i(t)$ and $\la_i(t)$ of $A(t)$ and $i = 1,...,n$  we obtain
\begin{equation}\label{errdeexpli}
  A(t) \dot x_i(t) - \la_i(t) \dot x_i(t) -  \dot \la_i(t)x_i(t) =    - \eta A(t)x_i(t) + \eta \la_i(t)x_i(t)  - \dot A(t) x_i(t) 
 \end{equation}
since scalars and vectors always commute, which was not the case for the $n$ by $n$ matrix equation (\ref{errdeexpl2}).\\
 Note that we do not know how to solve the full system eigenequation (\ref{errdeexpl2}) via ZNN  directly.\\
 Upon rearranging terms in (\ref{errdeexpli}) we finally have
\begin{equation}\label{errdeexplii}
 ( A(t) - \la_i(t) I_n) \dot x_i(t) - \dot \la_i(t)  x_i(t)  =   ( - \eta (A(t) - \la_i(t) I_n)   - \dot A(t)) x_i(t) 
 \end{equation}
where $I_n$ is the identity matrix of the same size $n$ by $n $ as $A(t)$.\\[1mm]
For each $i = 1,...,n$ the last equation (\ref{errdeexplii}) is a differential equation in the unknown eigenvector $x_i(t) \in \RR^n$ and the unknown eigenvalue $\la_i(t) \in \RR$ which  we rewrite   in augmented matrix form by further rearrangement as 
 \begin{equation}\label{errdeexplim}
\begin{pmat}
 A(t) - \la_i(t) I_n & -x_i(t)\\
 2 x_i^T(t) & 0 
 \end{pmat}
 \begin{pmat} \dot x_i(t)\\ \dot \la_i(t) \end{pmat}
    =   \begin{pmat} ( - \eta (A(t) - \la_i(t) I_n)   - \dot A(t)) x_i(t) \\ -\mu (x_i^T(t)x_i(t) -1) \end{pmat} \ .
 \end{equation}
 Here the second block row in (\ref{errdeexplim}) has been added below equation (\ref{errdeexplii})  by expanding the exponential decay differential equation for $e_2(t) = x_i^T(t) x_i(t) -1$ that is meant to ensure unit eigenvectors $x_i(t)$ throughout for a separate decay constant $\mu$.\\[1mm] Note that the leading $n$ by $n$ block matrix $A(t) - \la_i(t) I_n$ of the system  matrix in equation (\ref{errdeexplim}) is symmetric for symmetric matrix flows $A(t)$. To help speed up our computations we transform the error differential equation $2x_i(t)^T \dot x_i(t) = -\mu \cdot (x_i^T(t)x_i(t) -1)$ for $e_2$ slightly by dividing it by $-2$. That leads to the equivalent error function differential equation $-x(t)^T \dot x_i(t)  = \mu/2 \cdot (x_i^T(t)x_i(t) -1)$ for $e_2$. Thus for all symmetric matrix flows $A(t)$, we can  replace the non-symmetric ZNN model (\ref{errdeexplim}) by its symmetric version 
  \begin{equation}\label{errdeexplimsym}
\begin{pmat}
 A(t) - \la_i(t) I_n & -x_i(t)\\
 -x_i^T(t) & 0 
 \end{pmat}
 \begin{pmat} \dot x_i(t)\\ \dot \la_i(t) \end{pmat}
    =   \begin{pmat} ( - \eta (A(t) - \la_i(t) I_n)   - \dot A(t)) x_i(t) \\ {\mu}/{2} \cdot (x_i^T(t)x_i(t) -1) \end{pmat} \ .
 \end{equation}
  The full eigenvalue and eigenvector problem (\ref{EVequat}) of time-varying symmetric matrix flows $A(t) \in \RR^{n,n}$ can now be solved for each of its $n$ eigenpairs separately by using the matrix differential equation (\ref{errdeexplimsym}) with a symmetric system matrix for $i = 1,...,n$ in turn.\\[-3mm]
  
  In this paper we study and solve (\ref{EVequat}) with the help of convergent look-ahead difference schemes of varying truncation error orders for discrete, smooth and non-smooth symmetric time-varying matrix $A(t_k)$ inputs.

\section{Solving the Zhang Neural Network Error Equation via Look-ahead and Backward Discretization Formulas}

Our model (\ref{errdeexplim}) has recently been solved in real-time by using the standard ODE solver {\it ode15s} of MATLAB, see \cite{ZYLUH2018}. In \cite{ZYLUH2018},  model (\ref{errdeexplim}) was shown to be convergent and robust against data perturbations. %However, even with the ODE solver's relative error set to its smallest acceptable value of  $10^{-13}$ in MATLAB, the relative errors of the ODE based outputs for the eigenvalue equation $AV = VD$  never dipped below around $10^{-6}$ in our tests. Moreover, adaptive ODE solvers such as {\it ode15s} retrace their path  when large changes occur in the computed solution. Therefore general  ODE solvers seem not suitable for sensor activated data that is generated at fixed sampling times such as at  each 1/50th, 1/100th or 1/1000th of a second. This is so because to function properly, adaptive ODE solvers must interpolate and recompute the solution from data points at intermediate instances that may or may not be available from the sensors' output. 
%But in many real-world uses, system data  is only readily available at  predetermined discrete time instances. 
ODE based solvers can unfortunately not  be adapted easily to real-world sensor driven applications since they rely on intermediate data that may not be available with discrete-time sensor data. They work best for model testing of function valued time-varying $A(t)$ inputs. 
Thus there is a need to develop alternate discrete-time solvers such as  ZNN methods that go beyond what was first explored in continuous-time in \cite{ZYLUH2018}.\\[1mm]
For discretized symmetric input data $A(t_k)= A(t_k)^T \in \RR^{n,n}$ and $k = 0,1,2,...$ it is most natural to discretize the differential equation (\ref{errdeexplimsym}) in sync with the sampling gap $\tau = t_{k+1} - t_k$ which we assume to be  constant  for all $k$.   With choosing $\mu = \eta$ we have 
\begin{equation}\label{Pzqdef}
\begin{array}{c}
 P(t_k) = \begin{pmat}
 A(t_k) - \la_i(t_k) I_n & -x_i(t_k)\\
 - x_i^T(t_k) & 0 
 \end{pmat} \in \RR^{n+1,n+1} , \  \ \  z(t_k) = \begin{pmat}  x_i(t_k)\\  \la_i(t_k) \end{pmat} \in \RR^{n+1} \ , \\[3mm]
 \text{ and }  \ q(t_k) = \begin{pmat} ( - \eta (A(t_k) - \la_i(t_k) I_n)   - \dot A(t_k)) x_i(t_k) \\ {\eta}/{2} \cdot (x_i^T(t)x_i(t) -1)  \end{pmat}  \in \RR^{n+1}
 \end{array}
 \end{equation}
 our model (\ref{errdeexplimsym}) becomes  the set of matrix differential equations 
 \begin{equation}\label{Pzqde}
 P(t_k) \dot z(t_k) = q(t_k) 
 \end{equation}
for $k = 0,1,2,...$, each equidistant discrete time step $0 \leq t_k \leq t_f$, and  each eigenpair $x_i(t_k)$ and $\la_i(t_k)$ of $A(t_k)$. Note that $P(t_k)$ is always symmetric if the input matrix $A(t_k)$ is.\\[0mm]
 
{\bf \large The Discretization Process}\\[2mm]
The  differential equations (\ref{errdeexplimsym}) and (\ref{Pzqde}) are equivalent. They each  contain two derivatives: that of the unknown eigendata vector $z(t_k)$ and that of the input function $A(t_k)$.  In our discretization we replace the derivative $\dot A(t_k)$ by  the backward discretization formula
\begin{equation}\label{Aprimetk}
\dot A_k = \dfrac{11 A_k - 18 A_{k-1} + 9 A_{k-2} - 2 A_{k-3}}{6 \tau} 
\end{equation}
  of error order $O(\tau^3)$ \cite[p. 355]{EMU96} for example, where we have abbreviated $A(t_j)$ by writing $A_j$ for each $j$ and $\tau$ is the constant sampling gap from instance $t_j$ to $t_{j+1}$. Formula (\ref{Aprimetk}) is a four-instant backward difference formula. Our look-ahead finite difference formula for $z_k$ 
is the five-instant forward difference formula 5-IFD that we have adapted from \cite{LMUZb2018} for $\dot z_k$ as 
 \begin{equation}\label{5ifd}
\dot z_k = \dfrac{8z_{k+1}+z_k -6z_{k-1} - 5 z_{k-2} + 2 z_{k-3}}{18 \tau} \ .
\end{equation}
It also has the  truncation error order $O(\tau^3)$. 
By replacing $\dot z_k$ in (\ref{Pzqde}), multiplying equation (\ref{5ifd}) by $18\tau$ and reordering terms we obtain
\begin{equation}\label{Pdediscrete}
18\tau \cdot \dot z_k =8z_{k+1}+z_k -6z_{k-1} - 5 z_{k-2} + 2 z_{k-3} = 18 \tau (P \backslash q)
\end{equation}
where we have expressed the solution $x$ of the linear system $Px = q$ by the MATLAB symbol $ P \backslash  q$. By multiplying  equation (\ref{5ifd}) by  $18\tau$, we have raised  the truncation error order of the resulting difference equation (\ref{Pdediscrete})  by one power of $\tau$ to $O(\tau^4)$. Solving equation (\ref{Pdediscrete}) for $z_{k+1}$ gives us the following look-ahead finite difference equation for the discretized time-varying matrix  eigenvalue problem
 \begin{equation}\label{zkplus1}
 z_{k+1} = \dfrac{9}{4} \tau (P \backslash q) - \dfrac{1}{8} z_k + \dfrac{3}{4} z_{k-1} + \dfrac{5}{8} z_{k-2} - \dfrac{1}{4} z_{k-3}  \ .
 \end{equation}
Completely written out in terms of the block entries of $P$ and with $q$ as defined in  (\ref{Pzqdef}), the last equation becomes
\begin{equation}
\begin{array}{c}
z_{k+1} = -\dfrac{9}{4} \tau \left( \begin{pmat}
 A_k - \la_k I_n & -\hat z_k\\
 - \hat z_k^T & 0 
 \end{pmat}  \backslash 
 \begin{pmat} 
 \left( ( \eta (A_k - \la_k I_n) + \dfrac{11 A_k - 18 A_{k-1} + 9 A_{k-2} - 2 A_{k-3}}{6 \tau} \right)\hat z_k\\ 
 -\eta/2 \cdot(\hat z_k^T \hat z_k - 1) \end{pmat} \right) \hspace*{2mm} \\[6mm] \label{fulldiffequat}
  - \dfrac{1}{8} z_k + \dfrac{3}{4} z_{k-1} + \dfrac{5}{8} z_{k-2} - \dfrac{1}{4} z_{k-3}  \ .
  \end{array}
 \end{equation}
Here $\hat z_k = z_k^{(1,...,n)}= x_i(t_k) \in \RR^n$ consists of the leading $n$ entries of $z_k\in \RR^{n+1}$ and  contains  the current approximate unit eigenvector $x_i(t_k)$ for the $i$th eigenvalue $\la_i(t_k)$ of $A(t_k)$.  The five-instant forward difference formula 5-IFD (\ref{5ifd}) was discovered, developed 
and analyzed extensively in \cite[sections 2.3, 2.4]{LMUZb2018}. There is was shown to be zero stable and consistent and hence convergent when used in  multistep recursion schemes. \\
The same truncation error order as in (\ref{5ifd}) is achieved by the following convergent six-instant forward difference formula 6-IFD of \cite[Sections III A and Theorems 2 and 3]{QZY2019}
\begin{equation}\label{6ifd}
\dot z_k = \dfrac{13}{24\tau}z_{k+1} - \dfrac{1}{4\tau} z_k - \dfrac{1}{12\tau}z_{k-1} - \dfrac{1}{6\tau}
z_{k-2} - \dfrac{1}{8\tau}z_{k-3} + \dfrac{1}{12\tau} z_{k-4}
\end{equation}
with truncation error order $O(\tau^3)$. But it gives slightly better results when combined with the five-instance backward difference formula, see \cite[p. 356]{EMU96},
\begin{equation}\label{Aprimetk5}
\dot A_k = \dfrac{25A_k - 48 A_{k-1} + 36 A_{k-2} - 16 A_{k-3} + 3 A_{k-4}}{12\tau}
\end{equation}

\vspace*{-2mm}
\noindent  
of similar error order as our earlier discretization formula (\ref{Aprimetk}).\\[-3mm]
%Add 6-IFD here\\  %%%%%%%%%%%%%%%%%%%%

Very few convergent  look-ahead  finite difference formulas were  known until quite recently. The ones of higher  truncation error orders than Euler's method at $O(\tau^2)$  have been found by lucky and haphazard processes described  in \cite[section 2.3]{LMUZb2018} or \cite[Appendix A]{QZY2019} for example.  
Our five- and six-instant look-ahead difference formulas 5-IFD (\ref{5ifd}) and (\ref{6ifd}) with truncation error orders $O(\tau^4)$ serve rather well here with (\ref{6ifd})  and (\ref{Aprimetk5}) giving us one more accurate digit than  (\ref{5ifd})  and (\ref{Aprimetk}). 
In the meantime, higher truncation order look-ahead finite difference formulas have been developed \cite{U2019} with orders up to $O(\tau^8)$ and we will use some of these for comparisons in the next section with numerical examples.\vspace*{-4mm}

\section{Numerical Implementation, Results and Comparisons}

The ZNN algorithm that we detail in this section  finds the complete eigendata of  time-varying symmetric matrix flows $A(t) \in \RR^{n,n}$ by using the five-instance forward difference formula (\ref{5ifd})  and predicts the eigendata for time $t_{k+1}$ while using  the four instance backward formula (\ref{Aprimetk}) for approximating $\dot A(t_k)$. Using other difference  formula pairs such as  (\ref{6ifd})  and (\ref{Aprimetk5}) or higher error order ones instead  would  require finding additional starting values and adjusting the code lines that define $\dot A_k$ and $z_{k+1}$ accordingly.\\
To start the 5-IFD iteration  (\ref{fulldiffequat}) requires knowledge of  four start-up matrices $A(t_k) \in \RR^{n,n}$ for $ k =1,...,4$ and their complete  eigendata, i.e., knowledge of all   $n$  eigenvectors and associated eigenvalues of $A(t_{..})$ at  the  four time instances $t_1,...,t_4$. In this paper's  ZNN eigen algorithms we always  start from $t_1 = 0$  with $A(0)$ and then gather the complete eigendata of $A(0), A(t_2), A(t_3)$, and  $A(t_4)$  in an $n\!+\!1$ x 4 eigendata matrix {\tt ZZ} using Francis' QR algorithm as implemented by MATLAB's  {\tt eig}  m-file. To evaluate $z(t_{k+1}) \in \RR^{n+1}$ with an eigenvector preceding the respective eigenvalue for any $k \geq 4$  via formula (\ref{fulldiffequat}), we always rely on the four immediately time-preceding eigendata sets.\\[1mm]
 We repeat this ZNN process in a double do loop, in an outer loop for increasing time instances $t_k$ from $t_5$ on until the final time $t_f$ is reached, and in an the inner, a separate loop to  predict each of the $n$ eigenvalues and the associated eigenvectors  of $A(t_{k+1})$ from earlier eigendata  for the next time instance $t_{k+1}$ by using the look-ahead finite difference formula (\ref{5ifd}).\\[1mm]
This works well for smooth data changes in $A(t)$ and all eigenvalue distributions of the flow. Note that ZNN succeeds even for symmetric matrix flows with repeated eigenvalues, where both, decomposition and ODE path following methods suffer \emph{crucial breakdowns}, see e.g. \cite{DE99} and \cite{LM}. When the computed approximate derivative $\dot A(t_k)$ computed by a discretization formula such as (\ref{Aprimetk})  has an unusually large norm above 300, we judge that the input data $A(t_k)$ has undergone a discontinuity, such as from sensor or data transmission failure. In this case  we restart the process anew from time $t_k$ on, just as we did at the start from $t_1 = 0$ on  and again use Matlab's {\tt eig} function for four consecutive  instances to produce the necessary  start-up eigendata. Thereafter we  compute recursively again  from $t_{k+4}$  onwards via (\ref{fulldiffequat})  as done  before for $t_1 = 0, t_2,t_3$, and $t_4$, but now for the jump- modified input data flow $Aj(t)$. This proceeds until another data jump occurs and the original data flow $A(t)$ continues as input till the final time $t_f$. In that way we can find the complete eigendata  vectors $z_m(t_j)$ and $j > k+3$ for all affected indices $m$ until we reach the final specified time instance $t_f$. Our method thus allows for piecewise smooth data flows $A(t)$ and renders it adaptive to real-world occurrences and implementations.\\[1mm]
In the following we show some low dimensional test results for the symmetric time-varying  eigenvalue problem with the 7 by 7 seed matrices%\\[-4mm]

{\small
\vspace*{-5mm}
$$
A_s(t)  = \begin{pmat} 
 \sin(t)+2 &   e^{\sin(t)} & 0& -e^{\sin(t)}&           1/2&       1+\cos(t)  & 0\\
  e^{\sin(t)} &    \cos(t)-2&0&             1&        \cos(2t)&      1          & 0\\
            0    &           0&  -0.12t^2+2.4t-7& 0&   0&   0&          0\\
 -e^{\sin(t)}&                1&  0&    1/(t+1)&    \arctan(t)&    \sin(2t)&0\\
            1/2&       \cos(2t)&  0&\arctan(t)&                1&    e^{\cos(t)}&0\\
    1+\cos(t)&               1&    0&  \sin(2t)&    e^{\cos(t)}&   1/(t+2) &0\\
    0&0&0&0&0&0&      -0.15t^2+3t-6
         \end{pmat} 
  $$}
%   \vspace*{-8mm}
   
  \noindent
%  \vspace*{-8mm}
   and the data-jump perturbed symmetric time-varying  matrices \vspace*{-2mm} 
 
 {\small
  $$
A_{sj}(t) = \begin{pmat}  
      \sin(t)+2 &   e^{\sin(t)} & 0& -e^{\sin(t)}&        \tcr{0}&        1+\cos(t)&0\\
      e^{\sin(t)}&         \tcr{0}&  0&            1&         \cos(2t)&        1             &0\\
      0&                      0& \tcr{1.3t-15}&       0&              0&               0 &         0\\
     -e^{\sin(t)}&               1&  0&     1/(t+1)&          \tcr{1}&       \tcr{2\cos(2t)}&0\\
            \tcr{0}&      \cos(2t)&    0&     \tcr{1}&          \tcr{-3}&      e^{\cos(t)}&0\\
         1+\cos(t)& 1&0&             \tcr{2\cos(2t)}&   e^{\cos(t)}&      \tcr{6}/(t+2)&0\\
          0&0&0&0&0&0&           \tcr{14.05-t}
         \end{pmat}  \in \RR^{7,7}
    $$}
    
 \vspace*{-4mm}
 \noindent
whose entries we randomize via the same random orthogonal similarity $U \in \RR^{7,7}$  to become our test matrices $A(t) = U^TA_s(t)U$ and $Aj(t) = U^TA_{sj}(t)U$.
In our tests, $A_j(t)$  takes over from $A(t)$ between two-user set discontinuous or jump instances. $A_{sj}(t)$'s data differs from  that of $A_j(t)$ in positions (1,5), (2,2), (3,3), (4,5), (4,6), (5,1), (5,4), (5,5), (6,4), (6,6), and (7,7). These entries  are marked in red in $Aj(t)$.\\[1mm]
Our first example uses the sampling gap  $\tau$ = 1/200 sec  and the decay constant $\eta = 4.5$. The computations run from time $t$ = 0 sec to $t = t_f = 20$ sec with data input discontinuities at two jump points $t $ = 8 sec and $t$ = 14.5 sec. The computed eigendata satisfies the eigenvalue equation $A(t)V(t) = V(t)D(t)$ with relative errors in the $10^{-5}$ to $10^{-7}$ range in Figure 2 except for a sharp glitch near the end of the run when two eigenvalues of $A(t)$ take sharp turns and nearly cross paths.  Figure 3 shows the behavior of the norm of the derivative matrix $\dot A(t)$  as computed via (\ref{Aprimetk}) over time. This norm varies very little between 2 and 4 unless the input data becomes discontinuous at the chosen jump points where it spikes to around 2,000. Figure 4 examines the orthogonality of the computed eigenvectors over time. They generally deviate from exact orthogonality by between $10^{-4}$ and $10^{-7}$ except for spikes when two computed eigenvalues make sharp turns, sharper than the chosen sampling gap $\tau$ allows us to see clearly. The troublesome  behavior of 'coalescing' eigenvalues in matrix flows was noticed earlier by Mailybeav \cite{M06} and studied more recently by Dieci and Eirola \cite{DE99}. Without graphing, our algorithm's typical run times for this data set and duration until $t = t_f = 20$ sec   are around 0.5 seconds and take  about 2.5 \% of the total simulationtion time. This leaves the processor idle for around 97 \% of our 20 second time interval and indicates a very efficient real-time realization.\enlargethispage{50mm}

\vspace*{-38mm}
\begin{center}
\includegraphics[width=120mm]{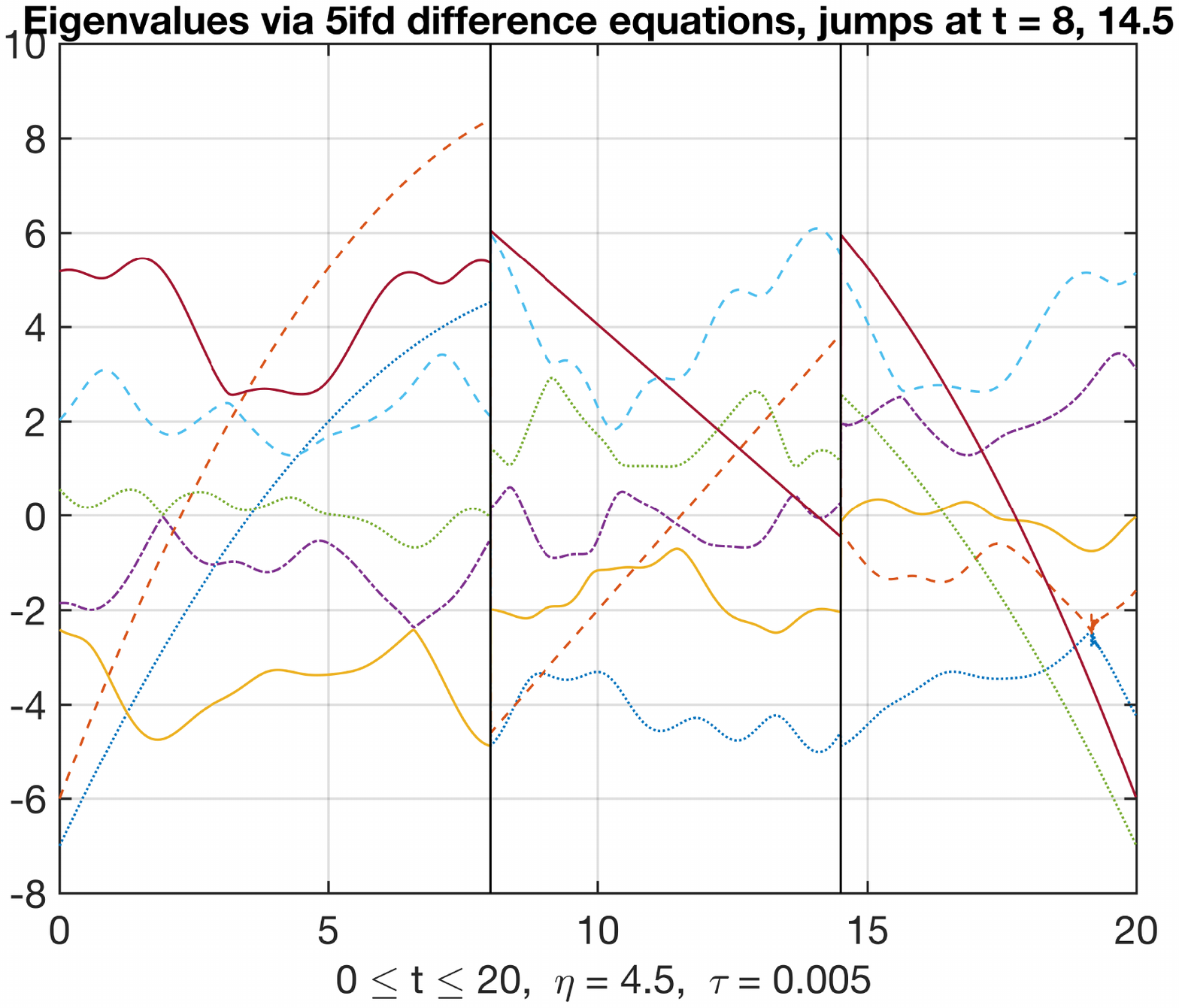}\\[-36mm]
Figure 1
\end{center} 
% \newpage
 
 %\enlargethispage{60mm}
 \vspace*{-41mm}   
  \begin{center}
\includegraphics[width=120mm]{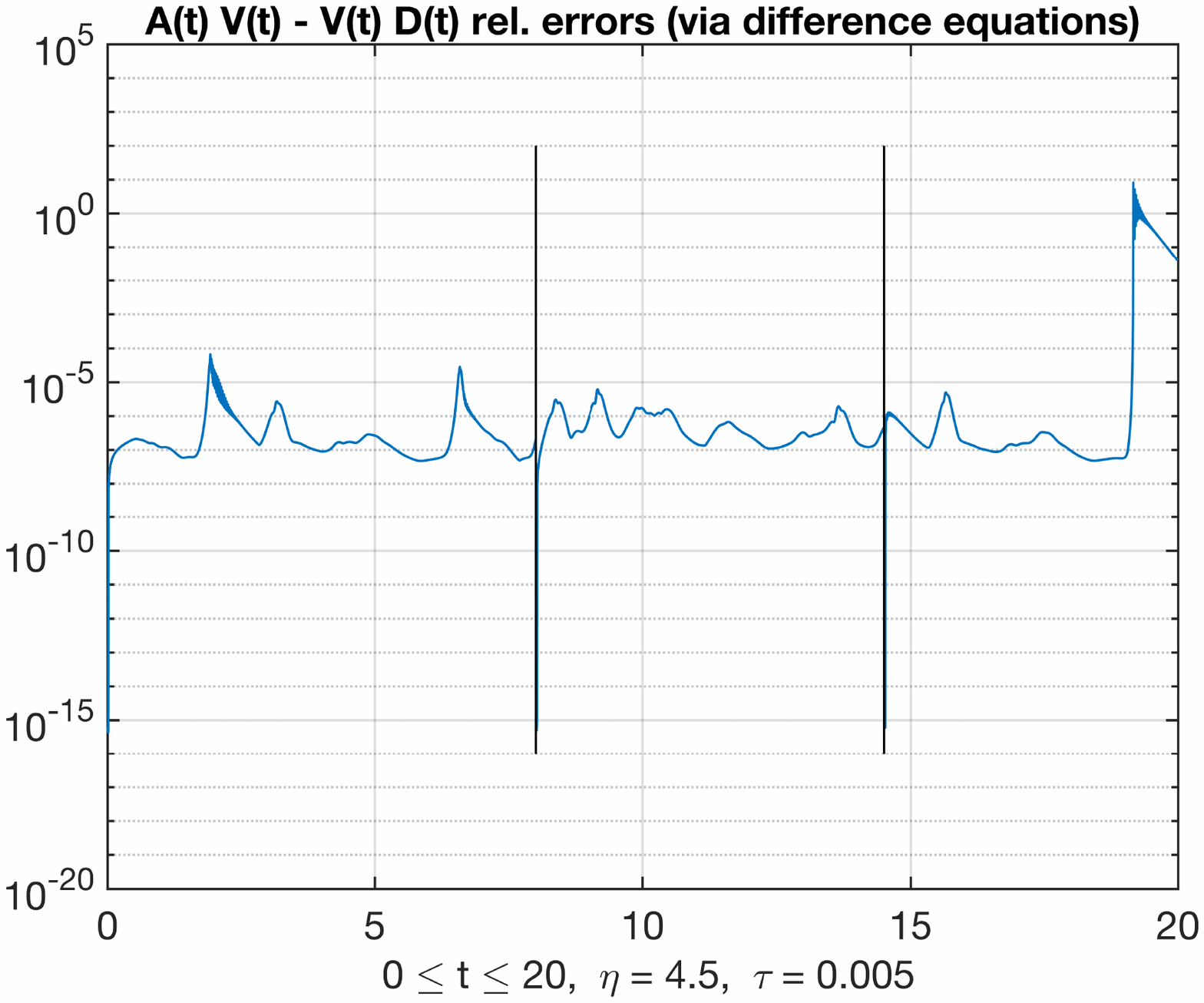}\\[-37mm]
Figure 2
\end{center}  
\newpage

\vspace*{-42mm}
\begin{center}
\includegraphics[width=120mm]{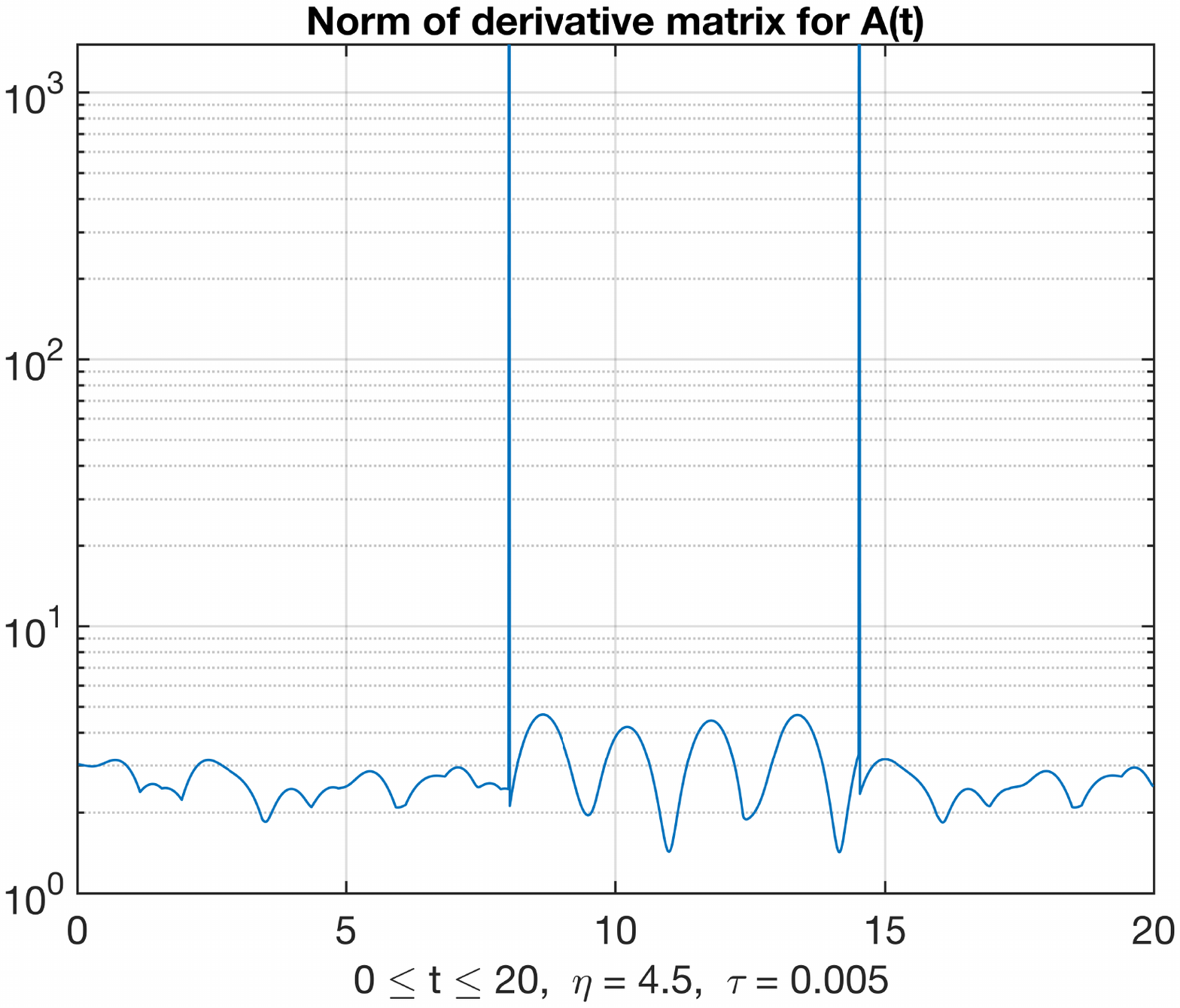}\\[-37mm]
Figure 3
\end{center}
%\newpage
 
\vspace*{-40mm}   
  \begin{center}
\includegraphics[width=120mm]{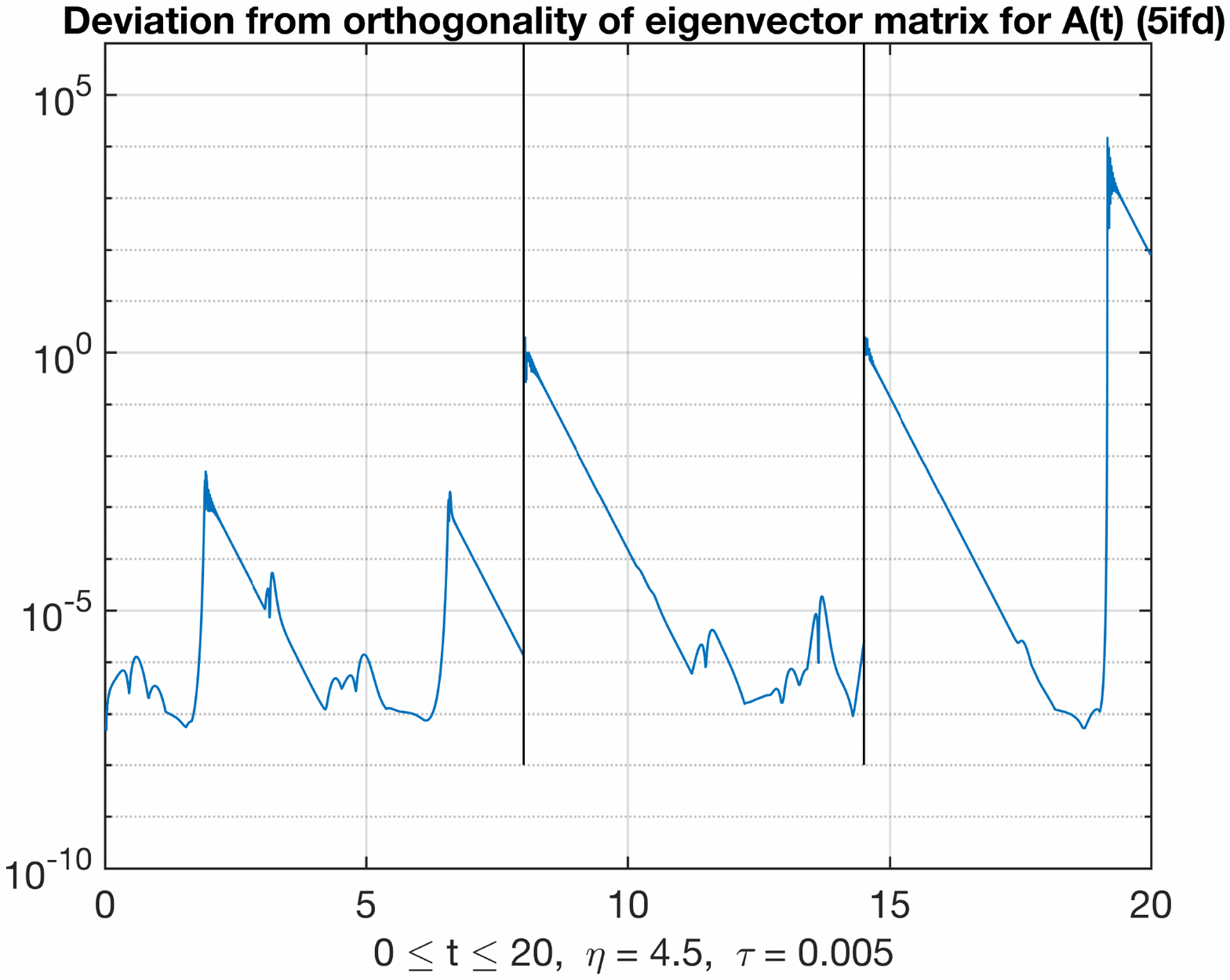}\\[-37mm]
Figure 4
\end{center}  
\newpage

\noindent
To reach higher accuracies with the sampling gap  shortened to $\tau$ = 1/1,000 sec, the results are much improved, see Figures 5 through 7 below. Here we have raised the value of $\eta$ to 80. When using its previous value of 4.5 with $\tau = 0.001$ sec, the plots would look much the same, but the achieved accuracy would suffer a wee bit.  On the other hand for $\tau = 0.005$ in Figures 1 -- 4, increasing $\eta$ from our chosen value  would not have worked well. With $\tau = 0.001$ sec and for any almost any $\eta\gg 1$, the average eigendata computation time for a run  without graphing  averages at around 2.5 sec. This allows the processor to be idle for around 87 \% of the total process time and again shows that our ZNN based discretized eigendata algorithm is well within real-time feasibility.
%\newpage

\enlargethispage{50mm}

\vspace*{-41mm}
\begin{center}
\includegraphics[width=120mm]{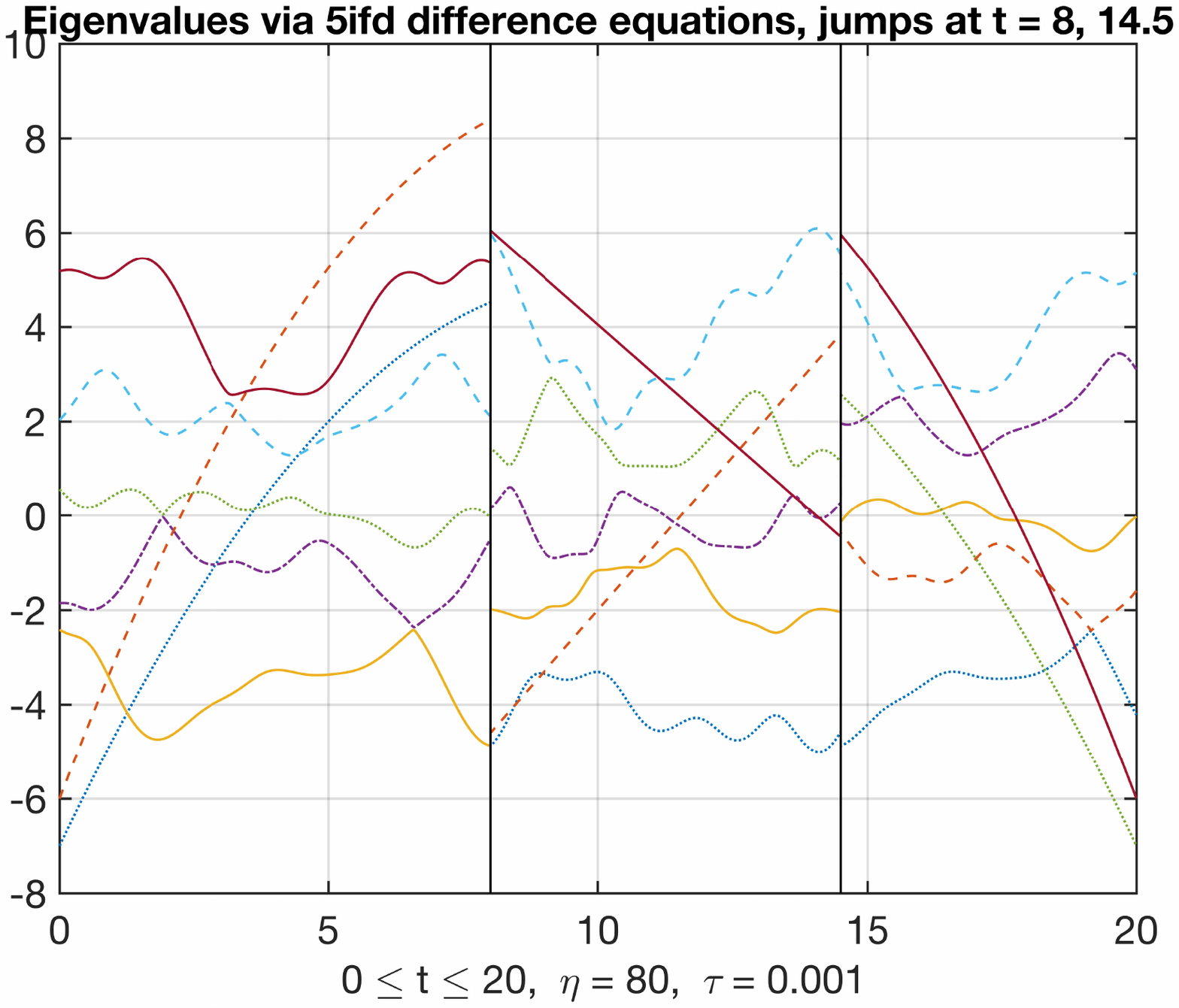}\\[-37mm]
Figure 5
\end{center} 

\vspace*{-2mm}
\noindent
Note that in Figure 5 the eigenvalue glitch of Figure 1 near --2 and $t = 19$ has been smoothed out. There are 31 incidences of eigenvalue crossings or repeated eigenvalues in this example that were handled by ZNN without any problems at all. \\[-4mm]

In Figure 6 below, the eigenvalue equation errors have been lowered by a factor of around $10^3$ by decreasing the sampling gap by a factor of 5 from Figure 2 and $5^4 = 625 \approx 1,000$, validating the 5-IFD error order $O(\tau^4)$ here.

 \vspace*{-43mm}   
  \begin{center}
\includegraphics[width=120mm]{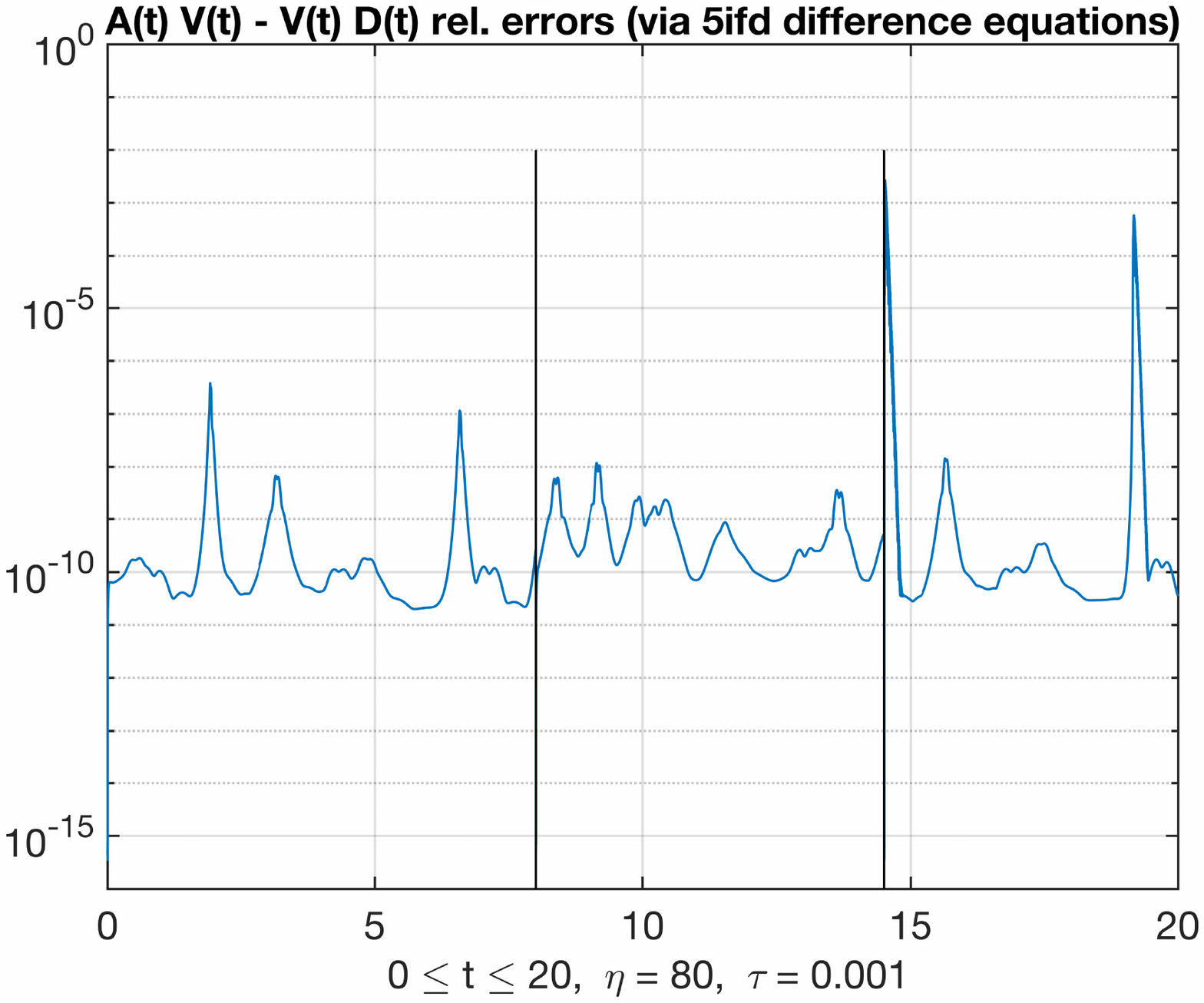}\\[-37mm]
Figure 6
\end{center}  
\newpage
 
% \enlargethispage{50mm}
\vspace*{-42mm}   
  \begin{center}
\includegraphics[width=120mm]{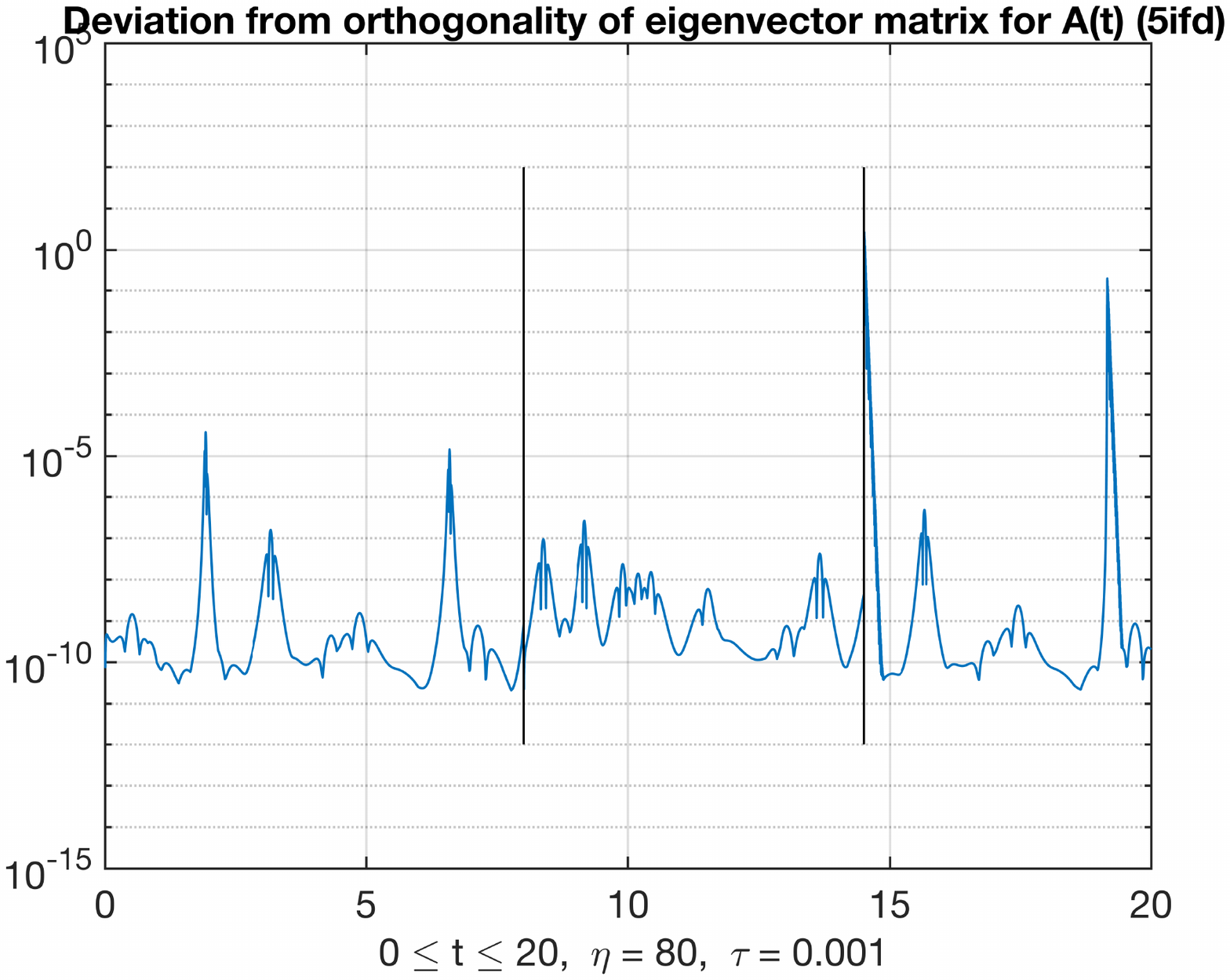}\\[-37mm]
Figure 7
\end{center}

The seven plots above were obtained by using the 5-IFD 1-step ahead difference formula (\ref{5ifd}) in conjunction with  the 4 instance backward difference formula (\ref{Aprimetk}). The 6-IFD (\ref{6ifd})  and the 5 instance backward formula (\ref{Aprimetk5}) pair have the same error order as the former difference formula pair (\ref{5ifd}) and (\ref{Aprimetk}) but they achieve superior results as the following plots indicate.
%\newpage

\enlargethispage{50mm}

\vspace*{-39mm}
\begin{center}
\includegraphics[width=120mm]{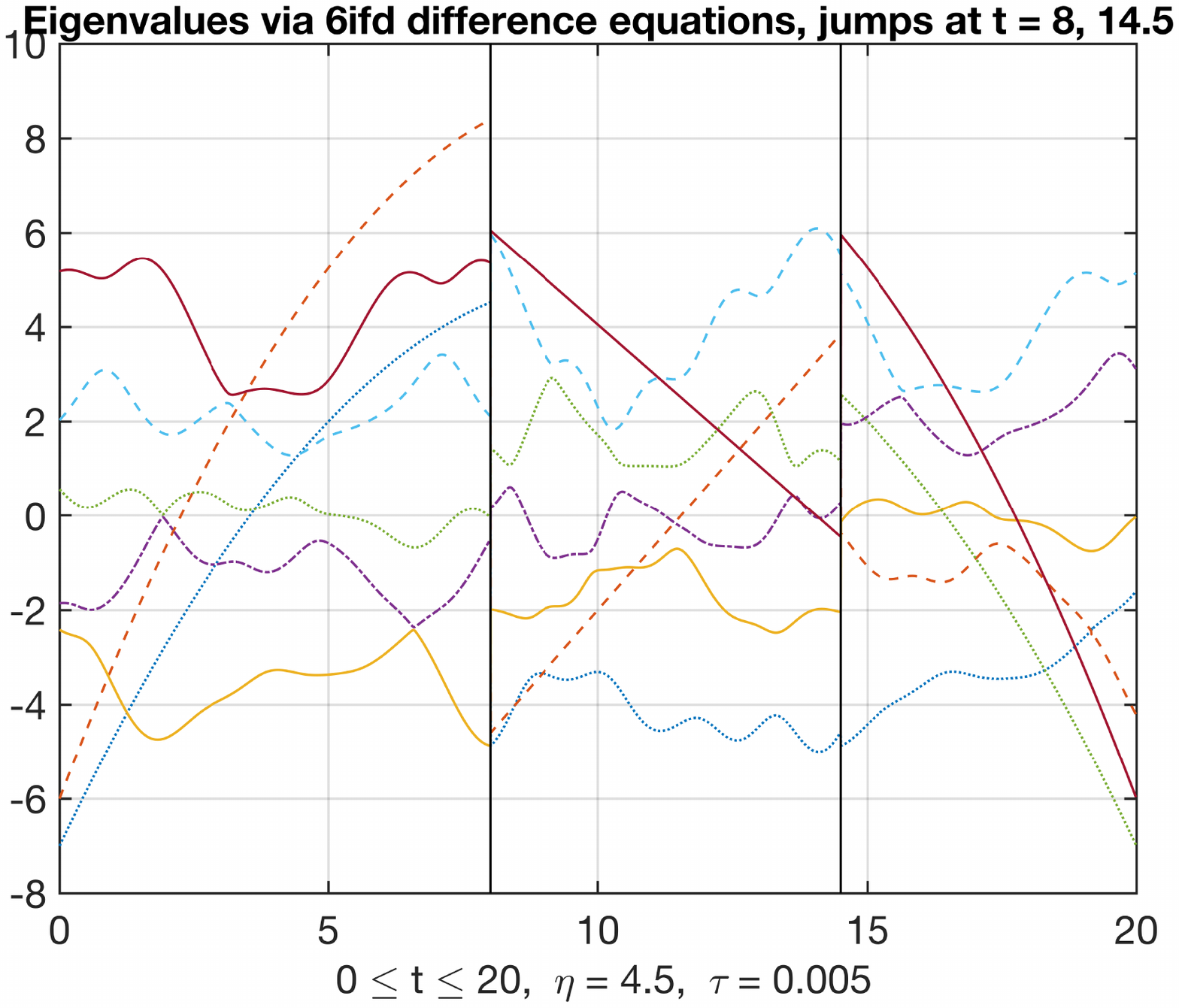}\\[-37mm]
Figure 8
\end{center}
This plot shows no glitches near --2 and $t = 19$ sec when compared with Figure 1 and it shows a clear 32nd  eigenvalue crossing at that former glitch point.\\
\newpage

 \vspace*{-45mm}   
  \begin{center}
\includegraphics[width=120mm]{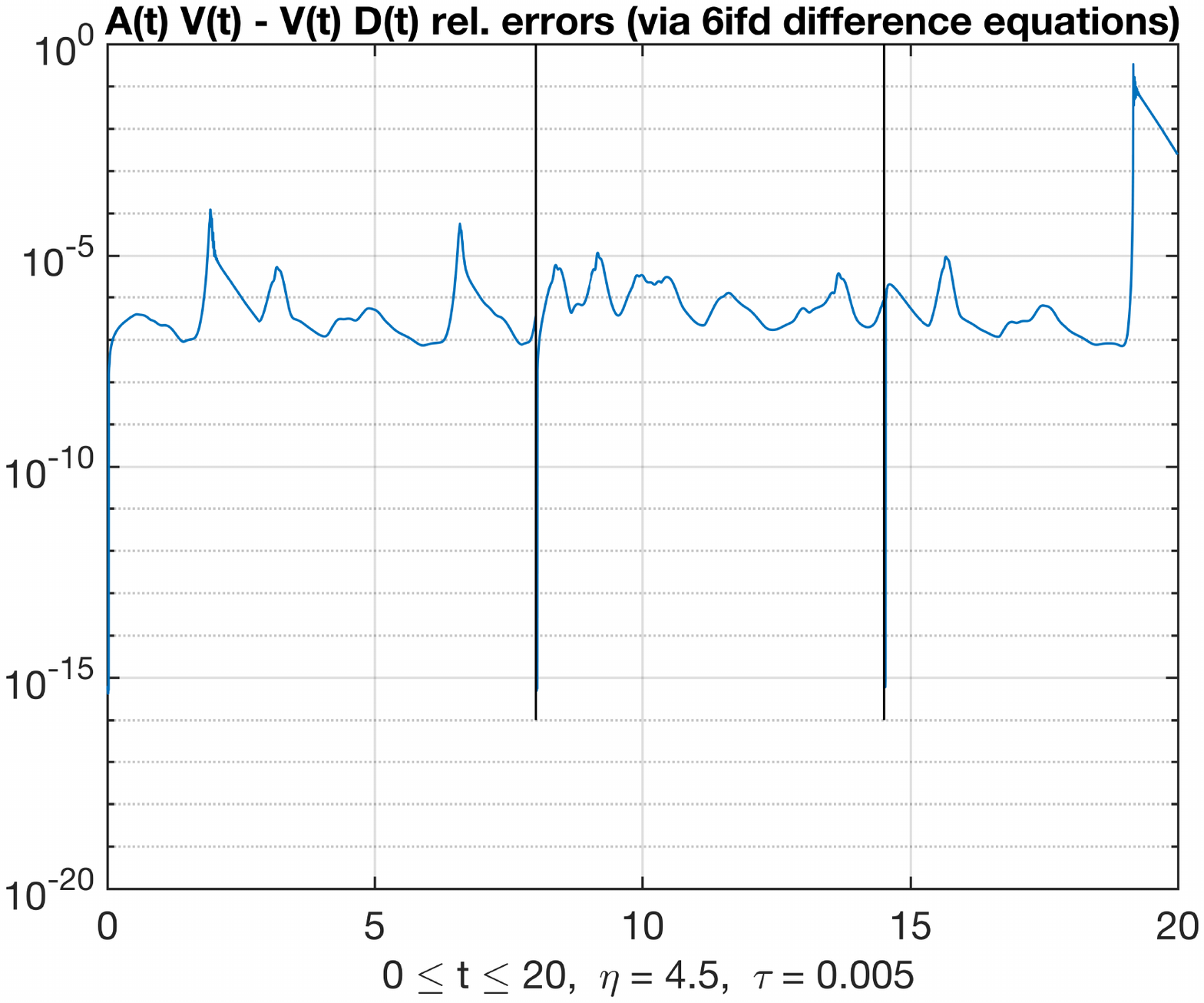}\\[-37mm]
Figure 9
\end{center} 

\vspace*{-3mm} \enlargethispage{60mm}
\noindent
 When using the 6-IFD in our Matlab code {\tt Zmatrixeig2\_3sym.m} \cite{U2018}, the computed results generally have lower relative errors of between a half to one digit than with the 5-IFD method and {\tt Zmatrixeig2\_2sym.m} in \cite{U2018}.    Both methods have the same truncation error order 4 and run equally fast.\\[-3mm]

%For further comparisons, we include two  eigendata graphs in Figures 10 and 11 below for the same test function $A(t)$ that are obtained by the earlier ODE solver based implementation of the  ZNN time-varying eigenvalue model, see \cite{ZYLUH2018}. However, as stated earlier we cannot incorporate $A(t)$ discontinuities into ODE based tests. The ZNN based {\tt ode15s} method computes the relevant eigendata in approximately half a second and this leaves the processor idle for about 97.5 \% of the time. Here the relative errors of the eigenvalue equations vary widely between $10^{-3}$ and $10^{-13}$, with most of the errors in the $10^{-7}$  range, as is depicted in Figure 11 below. This accuracy is much lower than the one depicted at around $10^{-10}$ in Figure 6 for the 5-IFD method for example.
%As mentioned before, the ODE based ZNN method \cite{ZYLUH2018} is not capable of running with constant sampling gap sensor inputs. In this test the ODE solver's largest integration step was $0.041$ sec and its shortest $0.000101$ sec with their average being around 0.0125 sec.
%\newpage
 
%\vspace*{-41mm}
%\begin{center}
%\includegraphics[width=120mm]{tveigode.pdf}\\[-37mm]
%\hspace*{2.6mm} Figure 10
%\end{center} 
%\newpage 

%\enlargethispage{60mm}
% \vspace*{-44mm}   
%  \begin{center}
%\includegraphics[width=120mm]{tvglrelerrode.pdf}\\[-37mm]
%\hspace*{2.6mm} Figure 11
%\end{center}  
%\newpage

Finally we investigate another way to try and predict the eigendata of $A_{k+1} = A(t_{k+1})$ from that of earlier eigendata for $A_j = A(t_j)$ with $j \leq k$. If the sampling gap $\tau$ is small such as $\tau = 0.001$ sec, how would the eigendata of $A(t_k)$ fare as a  predictor for that of $A(t_{k+1})$? To generate the plot below we have computed the eigenvalues of $A(t)$ 20,000 times using the static Francis  matrix eigenvalue algorithm that is built into Matlab in {\tt eig}.

\vspace*{-38mm}
\begin{center}
\includegraphics[width=110mm]{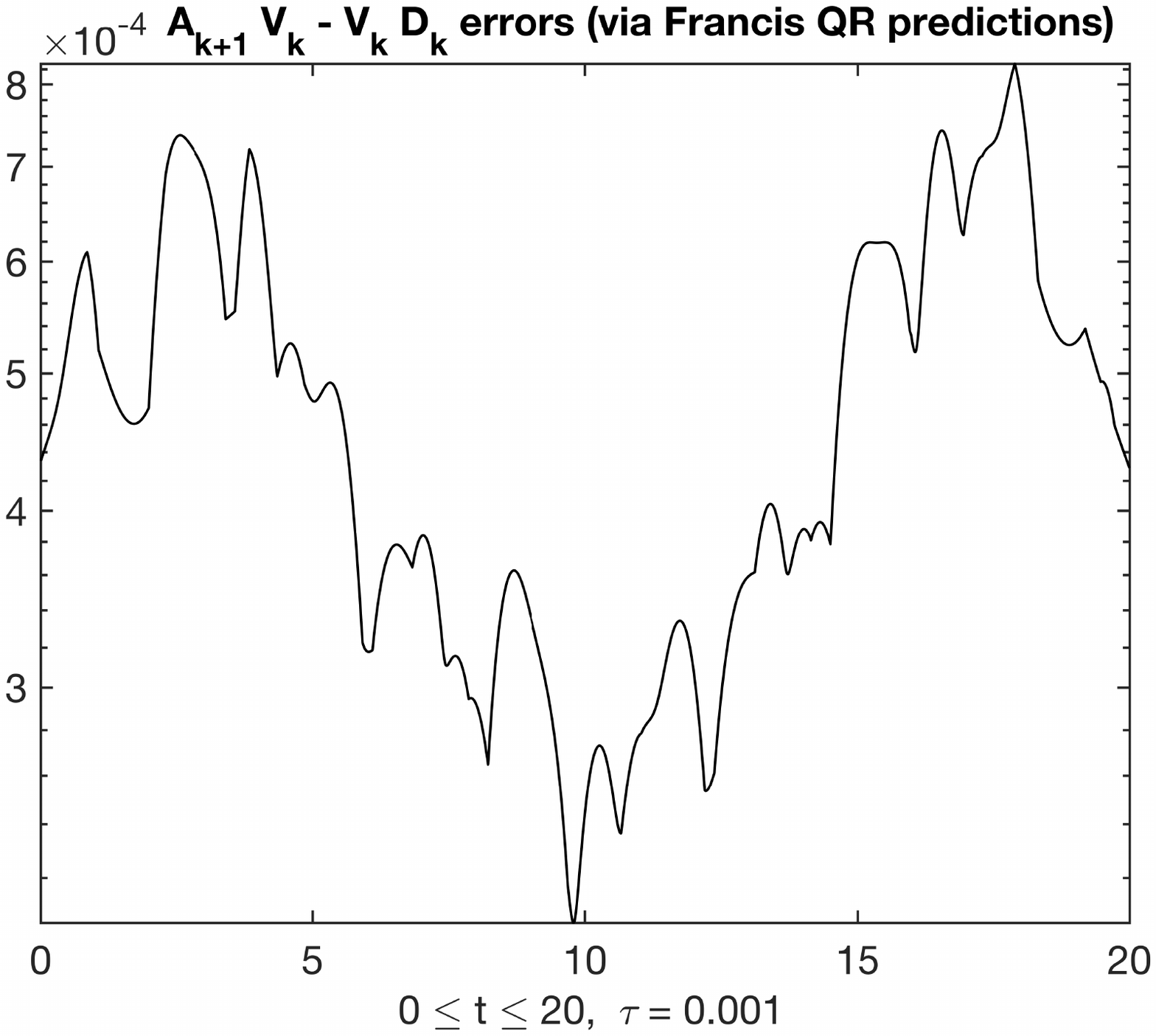}\\[-34mm]
\hspace*{2.6mm} Figure 10
\end{center} 

\vspace*{-2mm}
\noindent
Clearly this naive method does not reliably predict future eigendata at all since it generates past  eigendata whose average relative errors of magnitudes between $10^-4$ and  $8 \cdot 10^{-4}$ exceed those of our two tested ZNN algorithms significantly. This has been corroborated in general  for static methods in \cite{ZG2005}. Static method appropriations for time-varying problems suffer gfrom sizable
residual errors. The discretized ZNN methods of this paper are reliable predictors for time-varying matrix eigenvalue problems. ZNN can achieve  high accuracy. It is predictive and discretizable, as well as capable of handling discontinuous and non-smooth input data flows $A(t)$.
\newpage

\section*{Outlook}

In the course of this research we have come upon a number of interesting observations and intriguing open question which we post here.\\[2mm]
{\bf A :  \ \ Orthogonal Eigenvectors for Time-varying Symmetric Matrix Flows $A(t)$}\\[1mm] \enlargethispage{50mm}
Our ZNN based methods with low truncation error orders do not seem to compute nearly orthogonal eigenvectors for time-varying real symmetric matrix flows $A(t)$ as depicted  in Figures 4 and 7. In theory, the set of eigenvectors of every symmetric  matrix $A(t)$ can  be chosen as orthonormal  and the static Francis QR algorithm computes them as such by its very design of creating a number of orthogonal similarities on $A(t)$ that then converges to a diagonal eigenvalue matrix for the symmetric matrices $A(t)$.  But ZNN dynamics have no way of knowing that each $A(t)$ is symmetric and besides  they pay no heed to backward stable orthogonal static methods. We have tried to force orthogonality on the eigenvectors of each $A(t)$ by adjoining the following exponential decay differential equation as an additional  error function 
$$ee(t) = \begin{pmat} -& x_1 &-\\& \vdots&\\- & x_{k-1}& - \end{pmat} \begin{pmat} \vdots\\x_k(t) \\ \vdots \end{pmat} \stackrel{!}{=} o_{k-1}
$$
to the $P(t_k) \dot z(t_k) = q(t_k)$ formulation in equation (\ref{Pzqde}).  But this made  for worse orthogonalities than occur naturally with our simpler method. 
However, the new higher truncation error order convergent look-ahead  difference schemes  of \cite{U2019} seem to put this problem with discretized ZNN methods  to rest. See Figure 13 below where we have made use of the new 5th error order convergent 7-IFD finite difference method,  coded as {\tt Zmatrixeig3\_3bsym.m} in \cite{U2018} to give us
$$  z_{k+1} =  \dfrac{- 80z_{k} + 182z_{k-1} + 206z_{k-2} - z_{k-3} - 110z_{k-4} + 40z_{k-5}}{237} \ .
$$
\vspace*{-43mm}
\begin{center}
\includegraphics[width=120mm]{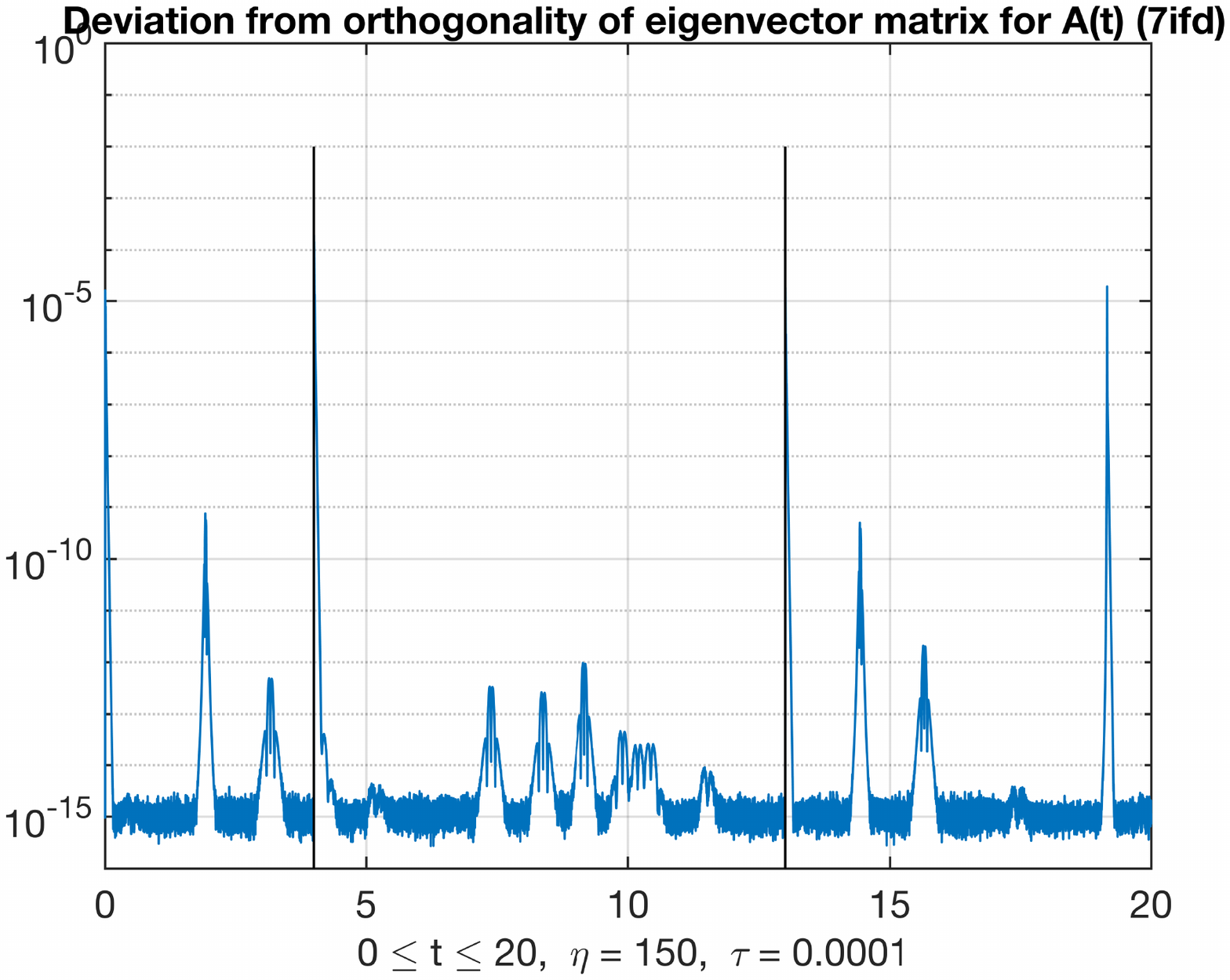}\\[-37mm]
\hspace*{2.6mm} Figure 11
\end{center} 
Figure 13 shows only moderate deviations from orthogonality of around $10^{-15}$ for the computed 7 by 7 eigenvector matrices  except for the deliberately chosen discontinuities  at 4.5 and 13. Otherwise there  only a couple of relatively small glitches with  orthogonality in the $10^{-9}$ to $10^{-12}$ range plus one major one in the $10^{-5}$ range near $t = 19$.
\newpage

{\bf B : \ \  General  Time-varying Matrix  Flows $A(t) \in \RR^{n,n}$ or $\CC^{n,n}$ and their Eigendata}\\[1mm]
How will the Zhang Neural Network method work for general square matrix flows $A(t)$, independent of their entry structure and eigenvalues? Its main task is updating the eigendata vector $z(t_k)$ by a  rule such as the  5-IFD based  formula below with truncation error order 4, coded as {\tt Zmatrixeig2\_2sym.m} in \cite{U2018}
$$ \hspace*{40mm}z_{k+1} = \dfrac{9}{4} \tau (P \backslash q) - \dfrac{1}{8} z_k + \dfrac{3}{4} z_{k-1} + \dfrac{5}{8} z_{k-2} - \dfrac{1}{4} z_{k-3} \ . \hspace*{30mm} (\ref{zkplus1})$$
This task is relatively simple:  create the matrix $P$ and the right hand side vector $q$ such as indicated  in formulas (\ref{Pzqdef}) and (\ref{fulldiffequat}), for example; compute a prescribed linear combination of previous $z$ vectors and solve a linear equation. By all signs, a generalization of our method should be able to work with arbitrary matrix flows, not just symmetric flows   $A(t) \in \RR^{n,n}$.  But we  have not  checked.\\[2mm]
\noindent
{\bf C : \ \ Matrix Flows $A(t) \in \RR^{n,n}$ for Large  Dimensions $n$, both Dense or Sparse, and their Eigendata\linebreak
\hspace*{8.5mm} Computations}\\[1mm]
How will an iteration rule such as  (\ref{zkplus1}) work with large dense or huge sparse matrix flows $A(t)$? Its main work consists of solving the linear system $P \backslash q$ for large dimensions $n$. Here all of our known static linear equation solvers for sparse and dense system matrices, be they Gauss, Cholesky, Conjugate Gradient, GMRES or more specific Krylov type methods should work well.  Future tests will tell. \\[2mm]
\noindent
{\bf D : \ \ Matrix Flow $A(t)$ Eigendata Computations with Singular System Matrices $P$ in Equation (\ref{Pzqde})\\
\hspace*{6mm}  $P(t_k) \dot z(t_k) = q(t_k)$ }\\[1mm]
Must ZNN based discrete methods break down if the DE system matrix $P(t)$ becomes singular at some time $t$? If $n$ is relatively small and we fail with Gaussian elimination to find $P \backslash q$,  i.e. if $P(t)$ is singular, then we can  append one row of $ n+1$ zeros to $P_{n+1,n+1}$  and a single zero to $q$,  and the backslash operator $\backslash$ of Matlab will not use elimination methods but solve the augmented non-square linear system 
$$ \hspace*{54mm} \begin{pmat}
& P & \\0& \dots & 0 
\end{pmat}_{n+2,n+1} \backslash \begin{pmat} q\\0\end{pmat}_{n+2} \hspace*{40mm} (\ref{Pzqde}\text{a})
$$ 
instead by using the SVD for the augmented equation (\ref{Pzqde}a). This is a little more expensive but delivers the least squares solution. We wonder what could be done for singular $P$ for sparse and huge time-varying matrix eigenvalue problems in the singular $P(t)$ case.\\[2mm]
{\bf E : \ \  Is it more advantageous to perform $n$ separate Iterations, each with one  $n+1$ by $n+1$ DE System Matrix\linebreak 
\hspace*{5.2mm}  $P$ as was done here, or use just one  $n^2 +n$ by $n^2 + n$ comprehensive DE System at each timestep $t_k$? }\\[1mm]
We have tried both approaches and the outcome depends. Further tests are warranted.\\[-2mm]

\hspace*{40mm} \underline{\hspace*{60mm}}\\[5mm]
\noindent 
Time-varying matrix  problems form a new and seemingly  different branch of Numerical Analysis. ZNN Neural Network based algorithms  rely on a few known stable 1-step ahead discretization formulas for their simplicity,  speed and accuracy, and on solving linear equations. ZNN method sensitivities such as the accuracy problems with sharply turning eigenvalue paths  as depicted in Figures 1, 2, 4, 5, 6, 7, 9 and 11 near $t = 19$ are new and differ from what we have learned in static matrix numerical analysis.  This gives us fertile ground for  in-depth numerical explorations of ZNN methods with time-varying matrix problems in the future. \\[-4mm]

\hspace*{40mm} \underline{\hspace*{60mm}}\\[3mm]
\noindent
Matlab codes of  our programs for  time-varying matrix eigenvalue problems   are available inside the folder \\    
\tcb{\tt  \url{http://www.auburn.edu/~uhligfd/m\_files/T-VMatrixEigenv}} \ . The 5-IFD version is  called  \ {\tt Zmatrixeig2\_2sym.m} (via 5-IFD),   the 6-IFD version is \ {\tt Zmatrixeig2\_3sym.m} (via 6-IFD) and the 7-IFD method is {\tt Zmatrixeig3\_3bsym.m}   in \cite{U2018}.\\[-1mm]

{\bf Acknowledgements : }  We thank Fran\c{c}oise Tisseur and Peter Benner and his group for valuable suggestions and additional references.%\\[-2mm]

 \vspace*{6mm}

\hspace*{30mm}{[} \ .../latex/t-vMatrixEvviaDiffEq.tex \ ]  \hspace{20mm}    \today

\vspace*{8mm}

\noindent
tveig005.pdf\\
tvglrelerr005.pdf\\
tvdernorm005.pdf\\
tvorthdev005.pdf\\[2mm]
tveig001.pdf\\
tvglrelerr001.pdf\\
tvorthdev001.pdf\\[2mm]
tv9eig005.pdf\\
tv9glrelerr005.pdf\\[2mm]
%tveigode.pdf\\
%tvglrelerrode.pdf\\[2mm]
tvglrelerrFrQR.pdf\\
tveig33symm.pdf

\end{document}